\documentclass[12pt]{article}
\usepackage{amscd, amssymb,latexsym,amsmath, amscd, amsmath}
\usepackage[all]{xy}
\begin{document}

\title{Relating invariant linear form and local epsilon factors
  via global methods} 
\date{}
\author{Dipendra Prasad}
\maketitle

\begin{abstract}
{We use the recent proof of Jacquet's conjecture due to Harris and Kudla,
and the Burger-Sarnak principle to give a proof about the
relationship between the existence of trilinear forms on representations
of $GL_2(k_u)$ for a non-Archimedean local field $k_u$ and local epsilon
factors which was earlier proved only in the odd residue characteristic
by this author in [P1]. The same method gives a global proof 
of a theorem of Saito and Tunnell about characters of
$GL_2$ using a theorem of Waldspurger about period integrals for $GL_2$ 
and also an extension
 of the theorem of Saito-Tunnell by this author in [P2] 
which was earlier proved only in odd
residue characteristic.}

\end{abstract}

\newtheorem{theorem}{Theorem}
\newtheorem{Proof}{Proof}
\newtheorem{lemma}{Lemma}
\newtheorem{prop}{Proposition}
\newtheorem{cor}{Corollary}
\newtheorem{example}{Example}
\newtheorem{conjecture}{Conjecture}
\newtheorem{definition}{Definition}
\newtheorem{question}{Question}
\newtheorem{memo}{Memo}
\newcommand{\rhobar}{\overline{\rho}}
\newcommand{\Sha}{{\rm III}}
\newtheorem{conj}{Conjecture}

{\sf 

\section{Triple products}
Let $\pi_1$, $\pi_2$ and $\pi_3$ be three irreducible admissible
infinite dimensional representations of $GL_2(k_u)$ for a non-Archimedean
local field $k_u$ with the product of their central characters trivial.
Let $D_u$ denote the unique quaternion division algebra over $k_u$.
For an irreducible admissible
discrete series representation $\pi$  of $GL_2(k_u)$, let $\pi'$
denote the representation of $D_u^*$ 
associated to $\pi$ by the Jacquet-Langlands correspondence,
and let $\pi'= 0$,  if $\pi$ is not a discrete series representation.

The author in his thesis [P1] studied the space of 
trilinear forms
$\ell: \pi_1 \otimes \pi_2 \otimes \pi_3 \rightarrow {\Bbb C}$ 
which are $GL_2(k_u)$-invariant.  Let $m(\pi_1 \otimes \pi_2 \otimes \pi_3)$
denote the dimension of the space of such trilinear forms, 
and let $m(\pi'_1 \otimes \pi'_2 \otimes \pi'_3)$ denote the
dimension of the space of $D_u^*$-invariant linear forms  
on  $\pi'_1 \otimes \pi'_2 \otimes \pi'_3$ (so 
$m(\pi'_1 \otimes \pi'_2 \otimes \pi'_3)$ is nonzero only if all the
$\pi_i$ are discrete series representations). 
We let $\epsilon(\pi_1 \otimes \pi_2 \otimes \pi_3) = 
\epsilon(\frac{1}{2}, \pi_1 \otimes \pi_2 \otimes \pi_3) $
denote the triple product epsilon factor defined by the Langlands-Shahidi 
method; under the condition that the
product of the central characters is  trivial, 
$\epsilon(\pi_1 \otimes \pi_2 \otimes \pi_3) = \pm 1$. Usually, the
epsilon factor depends on an auxiliary additive character of the field,
but in our case it is independent of it.

The main results
proved in [P1] were the following.

\begin{enumerate}

\item {\bf Multiplicity one  theorem:} $m(\pi_1 \otimes \pi_2 \otimes \pi_3) \leq 1$.

\item{\bf Dichotomy principle:} $m(\pi_1 \otimes \pi_2 \otimes \pi_3) 
+  m(\pi'_1 \otimes \pi'_2 \otimes \pi'_3) =
 1$.

\item {\bf Theorem about epsilon factors:} 
$m(\pi_1 \otimes \pi_2 \otimes \pi_3)= 1$ if and only if 
$\epsilon(\pi_1 \otimes \pi_2 \otimes \pi_3) =1$.

\end{enumerate}

A few words about the proofs. The multiplicity one theorem
was proved by the method of {\it Gelfand pairs}, as developed
by Gelfand and Kazhdan, and was thus based on general principles 
available to prove such theorems.

The dichotomy principle was eventually related to the character identity
$$\Theta_{\pi}(x) = -\Theta_{\pi'}(x),$$
at regular elliptic elements $x$ 
for discrete series representations $\pi$ and $\pi'$ of $GL_2(k_u)$
and $D_u^*$ respectively, associated to each other by the Jacquet-Langlands
correspondence; this combined with the theorem on finite groups
according to which 
$$  m(\pi'_1 \otimes \pi'_2 \otimes \pi'_3) = 
\frac{1}{{\rm vol}({D_u^*/k_u^*})} 
\int_{D_u^*/k_u^*}
\Theta_{\pi'_1}(g)\Theta_{\pi'_2}(g)\Theta_{\pi'_3}(g) dg,$$
and a suitable variant for $GL_2(k_u)$(!), proves the dichotomy
principle for supercuspidal representations, others being much more
straightforward by the orbit method of Mackey. 

These two theorems had reasonably satisfactory proofs. However, 
the theorem about epsilon factors was proved in a 
case-by-case way by reducing it 
to simpler epsilon factors studied by Tunnell [T], and this reduction
was possible only for representations of $GL_2(k_u)$ 
arising from characters of quadratic extensions of $k_u$, and was thus
incomplete in even residue characteristic. Also, the proof of the theorem
on epsilon factors given in [P1] left much to be desired as it was 
brute force calculation which essentially amounted to calculating
epsilon factor of the triple product and relating it to 
$m(\pi_1 \otimes \pi_2 \otimes \pi_3)$ through an explicit knowledge of the 
character of representations of $GL_2(k_u)$ which has been known for
a long time through the work of Sally and Shalika, the epsilon factor
(associated to Galois representations) 
being calculated  either directly, or through Tunnell's work  which also
was a brute force calculation with characters. (Later, there was the
elegant paper of Saito [Sa] which proved  Tunnell's theorem in all 
residue characteristics.)

The aim of this paper is to offer a global method for the result
on epsilon factors using the fact that the global analogue of the 
results on trilinear form, which is the period integral, 
has recently been proved by
 Harris and Kudla (this result  
was conjectured by H. Jacquet). More importantly, the present proof, unlike
the earlier one, offers a conceptual reason why the theorem on epsilon factors
holds good, and suggests that the natural proof
 of the general conjectures of Gross
and the author in [G-P] about local branching laws from $SO(n)$ to $SO(n-1)$ 
in terms of epsilon factors may consist of the following 3 steps.

\begin{enumerate}
\item Proof of the global conjecture in [G-P] about nonvanishing 
of the period integral in terms of $L$-value at $1/2$.

\item Proof of the conjecture about branching from $SO(n)$ to $SO(n-1)$
for unramified principal series representations of $SO(n)$.

\item Proof of the conjecture for ${\Bbb R}$.
\end{enumerate}

In the cases studied in this paper which consists of cases of $n \leq 4$,
all the 3 steps are known. However, for general $n$, {\it none} of the
3 steps are known at present.

We now recall the theorem of Harris and Kudla.

\begin{theorem}(Harris-Kudla)
Let $k$ be a number field, and 
$\Pi_1$, $\Pi_2$ and $\Pi_3$  three cuspidal automorphic 
representations of $GL_2({\Bbb A}_k)$ with the product of their
central characters trivial. For a quaternion division algebra D
over $k$, let $\Pi^D_i$ be the automorphic representations 
of $(D\otimes_k {\Bbb A}_k)^*$ associated to $\Pi_i$ by the
global Jacquet-Langlands correspondence, if it exists. Then
the central critical $L$-value 
$L(\frac{1}{2}, 
\Pi_1 \otimes \Pi_2 \otimes \Pi_3)$ is nonzero if and only if 
for some $D$, $\Pi^D_i$ exist as  automorphic representations
for which there are $f_i^D \in \Pi_i^D$ such that
$$\int_{D^* {\Bbb A}_k^*\backslash (D\otimes_k {\Bbb A}_k)^*} f_1^D(g) f_2^D(g) f_3^D(g) d^\times g \not = 0.$$
\end{theorem}

The proof of the above theorem does not use the theorem on epsilon factors,
so it is legitimate to use it to prove the following theorem.

\begin{theorem}

Let $\pi_1$, $\pi_2$ and $\pi_3$ be three irreducible admissible
infinite dimensional representations of $GL_2(k_u)$ for a non-Archimedean
local field $k_u$ with the product of their central characters trivial.
Then there exists a nonzero
$\ell: \pi_1 \otimes \pi_2 \otimes \pi_3 \rightarrow {\Bbb C}$ 
which is $GL_2(k_u)$-invariant
if and only if 
$\epsilon(\pi_1 \otimes \pi_2 \otimes \pi_3) =1$.

\end{theorem}

\noindent{\bf Proof}: If one of the representations $\pi_i$ is either 
principal series, or is a twist of the Steinberg representation,
then the result can be proved by simple calculations as done in [P1].
(Essentially because in all these cases except when one is dealing with
the triple product of the Steinberg, there is an invariant linear
form as follows by simple orbit methods; also, the epsilon factor
is easily calculated to be one in all these cases.) We will, therefore,
in the rest of the proof assume that all the representations
$\pi_i$ are supercuspidal representations of $GL_2(k_u)$ with
the product of their central characters trivial.

We will first prove that if   there is a 
nonzero 
$\ell: \pi_1 \otimes \pi_2 \otimes \pi_3 \rightarrow {\Bbb C}$ 
which is $GL_2(k_u)$-invariant, then 
$\epsilon(\pi_1 \otimes \pi_2 \otimes \pi_3) =1$. Then we
will use the dichotomy principle to say that if the space of
$GL_2(k_u)$-invariant forms on 
$\pi_1 \otimes \pi_2 \otimes \pi_3 $ 
is zero, then $\pi_i$ are all discrete series representations, and
for  the corresponding representations $\pi'_i$ of $D^*_u$, there
is a    $D^*_u$-invariant linear 
form on    $\pi'_1 \otimes \pi'_2 \otimes \pi'_3$. Now by the same
method employed to prove the first case (when there is a $GL_2(k_u)$-invariant
linear form), we prove that the local epsilon factor 
$\epsilon(\pi_1 \otimes \pi_2 \otimes \pi_3) $
in this case is $-1$. This will complete the proof of the theorem.

Fix a totally real number field $k$ and a 
place $u$ of $k$ such that the completion of $k$ at $u$ is the 
local field $k_u$ that we started with. 
Let $D$ be a quaternion division algebra over $k$ which is unramified 
at all the finite places, and for simplicity, we assume that $D$ is
ramified at all the infinite places. This can be done by choosing
$k$ appropriately so that it has even degree over ${\Bbb Q}$.
 
We let $\Pi_1$ and $\Pi_2$ be automorphic representations of
  $D({\Bbb A}_k)^*$ with local components $\pi_1$ and $\pi_2$ 
respectively at $u$, unramified at all the 
other finite places outside $u$. 
It is well known that local
supercuspidal representations can be obtained as the local component
of an automorphic representation which is unramified at all the other
finite places, and have some weights at infinity.

We will
choose $\Pi_3$ which will have $\pi_3$ as its local component at
$u$ (and some  infinity type). The
representation $\Pi_3$ will be so constructed that the period
integral 

$$\int_{D^*{\Bbb A}^*_k\backslash D^*({\Bbb A}_k)} f_1(g) f_3(g) f_3(g) d^\times g ,$$ 
is nonzero  for some choice of functions $f_i \in \Pi_i 
\subset L^2(D^*\backslash D^*({\Bbb A}_k))
$. We will state the
general lemma below, which is part of Burger-Sarnak philosophy,
which implies that such choices can be made.

Once we have automorphic representations $\Pi_i$ with nonvanishing
period integral, the theorem of Harris and Kudla implies that 

$$L(\frac{1}{2}, 
\Pi_1 \otimes \Pi_2 \otimes \Pi_3)\not = 0.$$
This implies, in particular, 
 that the global sign, $\epsilon( \frac{1}{2}, \Pi_1 \otimes \Pi_2 \otimes 
\Pi_3)$,  in the functional equation
for $L(s, \Pi_1 \otimes \Pi_2 \otimes \Pi_3)$ is 1. However,
the global sign in the functional equation is nothing
but the product of the local epsilon factors. 
From the information that $\Pi_1$ 
is a principal series representation  at all the finite places of $k$ 
except $u$ it is easy to see that the
epsilon factor at all the finite places except $u$ is 1: this just follows 
from the general fact that 
$\epsilon(\sigma) \epsilon (\sigma^\vee) = \det (\sigma)(-1)$.

The  triple product epsilon factor at all the infinite
places is $-1$ by invoking the corresponding theorem at infinity 
(the local representations at infinity are given to 
have invariant linear forms);
these calculations at infinity are simple consequences
of results about the decomposition of the 
tensor product of finite dimensional representations
of $D^*({\Bbb R})$, the so called Clebsch-Gordon theorem, 
and are  discussed in section 9 of [P1].

It is clear then from the factorisation
$$\epsilon( \frac{1}{2}, \Pi_1 \otimes \Pi_2 \otimes \Pi_3) 
= \prod_\mu 
\epsilon( \frac{1}{2}, \Pi_{1,\mu} \otimes \Pi_{2,\mu} \otimes \Pi_{3,\mu}) = 1,$$
and recalling that there are an even number of places at infinity 
that $\epsilon(\pi_1 \otimes \pi_2 \otimes \pi_3) =1$, thus proving that
if the space of $GL_2(k_u)$-invariant linear forms on 
$\pi_1 \otimes \pi_2 \otimes \pi_3$ is nonzero, the triple product epsilon 
factor  $\epsilon(\pi_1 \otimes \pi_2 \otimes \pi_3) $ is 1.

Now, by the Dichotomy theorem, 
when there is no $GL_2(k_u)$-invariant trilinear form, there will be
one on 
the corresponding representations $\pi'_i$ of $D^*_u$, and we can 
do exactly the same analysis by choosing a 
quaternion division algebra over $k$ which is ramified exactly 
at $u$, and at all the infinite places of $k$ which are now assumed by 
choosing the number field $k$ appropriately to be 
odd in number. Once again the sign in the functional equation
will be one, and since the theorem at infinity 
gives the sign $-1$ at each of the infinite places, we get in this
case $\epsilon(\pi_1 \otimes \pi_2 \otimes \pi_3) =-1$,
completing the proof of the theorem. Le us emphasize that 
in our proof the
local epsilon factor at a finite place is matched to one at infinity
by global means, reducing us to a much simpler problem.

\vspace{1cm}

\noindent{\bf Remark :} The epsilon factor used in [P1] were those arising
from Galois representations, whereas the epsilon factors used here are
those defined by the Langlands-Shahidi method. It is a theorem of 
Ramakrishnan, cf. theorem 4.4.1 of [Ra], that these two epsilon factors 
are the same.

\vspace{6mm} 

The following result is essentially due to Burger and Sarnak [BS], 
but not quite! A result
of this kind only for the infinite prime is Proposition 3.1 of [HL]. 
Adding finite primes causes no extra difficulty. However, for the
sake of completeness we give a self-contained proof of a result which is
adequate for our purposes, and which is totally transparent and elementary
consequence of the weak approximation theorem.

\begin{lemma}
Let $k$ be a number field, $S$ a finite set of {\it finite} 
places of $k$, $G$ a  reductive algebraic group defined
over $k$, and $H$ a reductive subgroup of $G$. Suppose that $Z$ 
is a central subgroup of $H$, which remains central in $G$ with the
property that $Z\backslash H$ has no nontrivial 
$k$-rational characters. Let 
$G_S = \prod_{v \in S}G(k_v)$, similarly let $H_S = \prod_{v \in S}H(k_v)$.
Let $\pi = \otimes_v \pi_v$ be an automorphic representation
of $G({\Bbb A}_k)$. Suppose that $\mu_v$ 
are supercuspidal representations of $H(k_v)$, $v \in S$, which are 
induced from representations $\nu_v$ of subgroups ${\cal K}_v$ 
which are certain open subgroups of $H_v$, compact modulo $Z_v$. 
Assume that $\mu_v$  
appears as a quotient of $\pi_v$ restricted to $H(k_v)$ for all $v \in S$. 
Then there
is an automorphic representation $\mu = \otimes_{v } \mu_v$ of $H({\Bbb A}_k)$
with $\mu_S = \otimes_{v \in S} \mu_v$, 
and  functions $f_1 \in \pi$,  $f_2 \in \mu$
such that 
$$\int_{H(k)Z({\Bbb A}_k) \backslash H({\Bbb A}_k)}f_1  
\bar{f}_2 dh \not = 0.$$

\end{lemma}

\noindent{\bf Proof :} By the assumption that $\mu_v$ 
appears as a quotient of $\pi_v$ restricted to $H(k_v)$ for all $v \in S$, 
the representation
$\nu_v$ of ${\cal K}_v$  is a subrepresentation of 
 $\pi_v$ restricted to ${\cal K}_v$ for all $v \in S$. 
This means that there
is a function $f$ on $G(k) \backslash G({\Bbb A}_k)$ whose ${\cal K}_S= 
\prod_{v \in S}{\cal K}_v$
translates generate a space of functions which is isomorphic to
$\otimes \nu_v$ as ${\cal K}_S$-modules. We prove that the restriction
of {\it one} such  function to $H(k) \backslash H({\Bbb A}_k)$ is not 
identically zero.
Observe that $G({\Bbb A}^S) = \{ x \in G({\Bbb A}), 
x = \prod_v x_v| 
x_v = 1, \forall v \in S \}$ operates on such functions (by right translation), and if all the
$G({\Bbb A}^S)$ translates of a function $f$ were zero at the identity 
element of $G({\Bbb A})$,  the function $f$ would be identically zero 
by the weak approximation theorem 
according to which $G(k)$ is dense in $G_S$.
Thus we have a function $f$ 
 on $G(k) \backslash G({\Bbb A}_k)$ 
whose restriction to $H(k) \backslash H({\Bbb A}_k)$, say $\tilde{f}$,
 is not zero. The ${\cal K}_S$
translates of $\tilde{f}$ generates a space of functions 
now on $H(k) \backslash H({\Bbb A}_k)$ which is isomorphic to
$\otimes \nu_v$ as ${\cal K}_S$-modules. Since,
$$\mu_S = {\rm ind}_{{\cal K}_S}^{H_S} \nu_S, $$ 
the $H_S$-translates of $\tilde{f}$  generates a space of functions 
on $H(k) \backslash H({\Bbb A}_k)$ which is isomorphic to
$\mu_S$ as ${H}_S$-modules.
We are now done by the next lemma.

\begin{lemma}
Suppose that $H$ is a reductive  algebraic group over a number 
field $k$ with center $Z$. Suppose
that $S$ is a finite set of finite places and 
$H_S = \prod_{v \in S}H(k_v)$. Suppose that $\tilde{f}$ is a
continuous function on  $H(k) \backslash H({\Bbb A}_k)$ 
with unitary  central character 
$\chi$ whose $H_S$ translates
generate an {\it irreducible} $H_S$-submodule, say $\mu_S$, 
of a  space of functions on  $H(k) \backslash H({\Bbb A}_k)$. 
Then there exists an automorphic function $g$ on $H({\Bbb A}_k)$ 
with unitary  central character  $\chi$ generating
an irreducible representation  with $H_S$-type $\mu_S$, and with 
$$\int_{H(k) Z({\Bbb A}_k)\backslash H({\Bbb A}_k)}\tilde{f}\bar{g} dh \not = 0.$$

\end{lemma}

\noindent{\bf Proof :} The proof of the lemma follows from the fact that the
space of automorphic forms on $H({\Bbb A}_k)$ is dense in the space of
continuous functions. The nonvanishing of the integral implies that the
$H_S$-type of the space generated by $g$ is the same as that of the space
generated by $\tilde{f}$ which is $\mu_S$.

\vspace{3mm}

\noindent{\bf Remark :} In the proof of theorem 2, we apply lemma 1
to G =  $GL_2 \times GL_2$, $H = GL_2$, and $Z = G_m$, using the
well-known theorem of Kutzko that a supercuspidal representation
of $GL_2(k_u)$ can be obtained as an induced representation from a finite
dimensional representation of an open subgroup which is compact modulo
center.

\vspace{6mm}

\section{Saito-Tunnell}

By exactly the same method as employed in the previous section, one 
can deduce the theorem of Saito and Tunnell which 
describes which characters of $L^*$, for $L$ a quadratic extension of a local
field $k_u$, appear in an  irreducible admissible
infinite dimensional representations of $GL_2(k_u)$, or in an 
irreducible representation of $D_u^*$, where $D_u$ is the unique
quaternion division algebra over $k_u$, in terms of the local epsilon factors.
It is elementary to see that characters of $L^*$ appear with multiplicity
at most 1 in any irreducible representation of $GL_2(k_u)$, or of $D_u^*$;
see for example, remark 3.5 of [P1].  The dichotomy principle too holds
in this situation, i.e., if $\pi$ is a discrete series representation 
of $GL_2(k)$, and $\pi'$ is the corresponding finite dimensional irreducible
representation of $D_u^*$, then for a character $\chi$ of $L^*$ whose
restriction to $k_u^*$ is the same as the central character of $\pi$, $\chi$ 
appears in exactly one of the representations $\pi$ or $\pi'$. This
follows from the character identity $$\Theta_{\pi}(x) = -\Theta_{\pi'}(x).$$

By appealing to a global theorem  due to
Waldspurger, we give a proof of the Saito-Tunnell theorem. 
However, there seems little  point in giving many details
except to recall the statements of the theorems.

\begin{theorem}(Waldspurger) Let $k$ be a number field, 
$F$ a quadratic extension
of $k$, $D$ a quaternion division algebra over $k$ containing $F$, 
and $\pi'$ an automorphic form on $D({\Bbb A})^*$ realised on a 
space of functions $E'$ on $D^*\backslash D({\Bbb A}_k)^*$. Let $T$ be the
torus inside $D^*$ defined by $F$, and $\Omega = \otimes \Omega_v$ a continuous
character of $T(k)\backslash T({\Bbb A}_k)$ whose restriction to 
the center of $D({\Bbb A})^*$ is the same as the central character
of $\pi'$. Then there exists a function $e'$ in $E'$ such that the
integral 
$$\int_{T(k){\Bbb A}_k^*\backslash T({\Bbb A}_k)} e'(t) \Omega^{-1}(t) dt,$$
is nonzero if and only if the following two conditions are verified:

\begin{enumerate}
\item For all places $v$ of $k$, the local representation $\pi'_v$
has $\Omega_v$-invariant linear form for the torus $(F \otimes k_v)^*$.
\item If $\Pi$ denotes the base change of $\pi$ to $GL_2({\Bbb A}_F)$,
$$L(\frac{1}{2}, \Pi \otimes \Omega^{-1}) \not = 0.$$
\end{enumerate}

\end{theorem}

We now state the theorem of Saito and Tunnell which follows from
the theorem of Waldspurger just as our earlier proof. In this theorem,
the epsilon factor is sensitive to the additive character chosen. 
The theorem is deduced from the theorem of Waldspurger which has
no reference to the additive character. The reason for this of course
is that the corresponding statement at infinity also depends on the choice
of such an additive character. (Note that fixing a character
of ${\Bbb A}_k/k$  at one place of $k$ fixes it also at any other place 
of $k$ because $k_v$ is dense in ${\Bbb A}_k/k$ for any place $v$ of $k$.) 
The theorem of  Saito and Tunnell although not stated in their
papers for  Archimedean fields is valid for such fields too, and has
quite an elementary proof. 

\begin{theorem}(Saito, Tunnell)
Let $\pi$ be an infinite dimensional irreducible admissible
representation of $GL_2(k_u)$, $L$ a quadratic extension
of $k_u$, and $\Pi$ the base change lift of $\pi$ to $GL_2(L)$.
Fix a nontrivial additive character $\psi$ of $L$ which restricted to 
$k$ is trivial.
Then for a character $\Omega_u$ of $L^*$ which has the same
restriction to $k_u^*$ as the central character of $\pi$,
$\epsilon(\Pi \otimes \Omega_u^{-1}, \psi)$ is indepenent of $\psi$ (as
long as  its restriction to $k_u$ is trivial), and takes the
value $\pm 1$.
The character $\Omega_u$ of $L^*$ appears in $\pi$ if and only if 
$$\epsilon(\Pi \otimes \Omega_u^{-1}, \psi)= 1.$$

\end{theorem}  

\section{Extending Saito-Tunnell}

Let $L$ be a  quadratic extension of a local field $k_u$
of characteristic $\not = 2$. The theorem of Saito and Tunnell
discussed in the previous section describes which characters of $L^*$
appear in an  irreducible admissible infinite dimensional 
representations of $GL_2(k_u)$, or in an 
irreducible representation of $D_u^*$ where $D_u$ is the unique
quaternion division algebra over $k_u$ in terms of the local epsilon factors.
If the representation $\pi_\theta$ of $GL_2(k_u)$ comes from a character
$\theta$ of $L^*$ via the construction of the Weil representation,
cf. [J-L, theorem 4.6], the representation $\pi_\theta$
decomposes into two irreducible representations $\pi_\theta= \pi_+ \oplus
 \pi_-$
when restricted to $GL_2(k_u)^+ = \{x \in GL_2(k_u)| \det(x) \in NL^* \}$
where $NL^*$ is the subgroup of $k_u^*$ of index 2 consisting of norms
from $L^*$; similarly one can define a subgroup of index 2 inside
$D_u^*$, to be denoted by $D_u^{*+}$, and for which one has a 
similar decomposition $\pi'_\theta= \pi'_+ \oplus \pi'_-$ of the corresponding
representation of $D_u^*$. 
Clearly $L^*$ is contained in both $GL_2(k_u)^+$, and $D_u^{*+}$, and one
can ask about a  generalisation of the theorem of Saito and Tunnell
to describe the decomposition of $ \pi_+, \pi_-, \pi'_+, \pi'_-$
restricted to $L^*$. Here is such a theorem.

\begin{theorem}
Let $\pi_\theta$ (resp. $\pi_{\theta}'$) be the irreducible admissible 
representation of $GL_2(k_u)$ (resp. $D_u^*$)
associated to a character $\theta$ of $L^*$. Fix embeddings of $L^*$
in $GL_2(k_u)^+$, and $D_u^{*+}$ (in  general there are two conjugacy classes
of such embeddings). Let $\psi$ be a nontrivial 
 additive character of $L$
trivial on $k_u$. The restriction of $\pi_\theta$ to  $GL_2(k_u)^+$ 
can be written as  $\pi_\theta=  \pi_+ \oplus \pi_-$ 
and the restriction of $\pi'_\theta$ to 
$D_u^{*+}$ can be written 
as $\pi'_\theta= \pi'_+ \oplus \pi'_-$, such that
 a character $\chi$ 
of $L^*$ with $(\chi \cdot \theta^{-1})|_{k_u^*} =\omega_{L/k_u}$
appears in $\pi_+$ if and only if $\epsilon(\theta \chi^{-1},\psi) 
=\epsilon(\bar{\theta} \chi^{-1},\psi) =1$, and appears in 
$\pi_-$ if and only if  $\epsilon(\theta \chi^{-1},\psi) =
\epsilon(\bar{\theta} \chi^{-1},\psi) =-1$. 
Similarly, a character 
 $\chi$ 
of $L^*$ with $(\chi \cdot \theta^{-1})|_{k_u^*} =\omega_{L/k_u}$
appears in $\pi'_+$ if and only if $\epsilon(\theta \chi^{-1},\psi) =1$,
and  $\epsilon(\bar{\theta} \chi^{-1},\psi) =-1$, and appears in 
$\pi'_-$ if and only if  $\epsilon(\theta \chi^{-1},\psi) = -1$, 
and $\epsilon(\bar{\theta} \chi^{-1},\psi) =1$. 

\end{theorem}

This theorem was proved in the odd residue 
characteristic case by this author in [P2]. We now prove this theorem 
in general by a global
argument similar to the one in the earlier sections. The following lemma
will play an important role in transferring information from a finite
prime of a number field to an infinite prime.

\begin{lemma} 
Let $k$ be a CM number field, and $u$ a place of $k$. Let $\theta:
k^*_{u} \rightarrow {\Bbb C}^*$ be a character. Then there exists a 
Gr\"ossencharacter $\Theta: {\Bbb A}_k^*/k^* \rightarrow {\Bbb C}^*$
which is unramified at all the finite primes outside $u$, has
$\theta$ as the local component at $u$, and whose local component
at infinity is given by 
$$\Theta_\infty(z_1,\cdots, z_r) = z_1^{n_1}\cdots z_k^{n_r}, ~~~~~{\rm for}
~~~ |z_i| = 1,$$
where $(n_1,\cdots, n_r)$ is any element of a certain subgroup of finite
index in ${\Bbb Z}^r$.
\end{lemma}

\noindent{\bf Proof :} Let $U_k= \prod_{v \in S_f} U_v \prod_{v\in S_\infty}
{\Bbb S}^1$ be the maximal compact subgroup of ${\Bbb A}_k^*$. 
Define a character
$\mu$ on $U_k$ by declaring its value on $U_{u}$ to be $\theta$ 
restricted to $U_{u}$, trivial on $U_v$, $v \not = u$, and given on
$(z_1,\cdots, z_r) \in ({\Bbb S}^1)^r$ to be $z_1^{n_1}\cdots z_r^{n_r}$.
Since $k^* \cap U_k = \mu_k$, the group of roots of unity in $k^*$,
$\mu|_{k^* \cap U_k}$ is a character of finite order which can be assumed to be
trivial by choosing the $r$-tuple of integers 
$(n_1,\cdots, n_r) \in {\Bbb Z}^r$ from a subgroup of finite
index in ${\Bbb Z}^r$. 

Thus we have a character, say $\Theta'$,
 of the group $U_k/ (k^* \cap U_k)$.
Since, $U_k/ (k^* \cap U_k) \hookrightarrow {\Bbb A}_k^*/k^*$,
and $U_k/ (k^* \cap U_k)$ is a compact group, in particular closed, 
the character $\Theta'$ 
can be extended 
to a Gr\"ossencharacter, say $\Theta''$, 
 of ${\Bbb A}_k^*/k^*$. Observe that at the moment
we have constructed a Gr\"ossencharacter $\Theta''$ 
of ${\Bbb A}_k^*/k^*$ which on 
$U_{u}$ is  $\theta$  restricted to $U_{u}$. However, we still have the
character $|\cdot |: {\Bbb A}_k^*/k^* \rightarrow {\Bbb C}^*$ given by
$x \rightarrow |x|$, which is trivial on $U_k$, and takes the value 
$q_{u}^{-1}$
on $\pi_{u} \in k_{u}^* \subset {\Bbb A}_k^*/k^*$. Therefore for an 
appropriate choice of $s_0 \in {\Bbb C}$, we can assume that
the Gr\"ossencharacter
$\Theta = \Theta''|\cdot |^{s_0}$ has $\theta$ as its local component at
$u$, is unramified at all the other finite places, and has 
the desired behaviour at infinity.

\vspace{1cm}

\noindent{\bf Proof of theorem 5:} Using the theorem of 
Saito-Tunnell, we note that this theorem is equivalent to
proving that
for a character $\chi$ of $L^*$ appears in an irreducible component of
 $\pi_\theta$ restricted to $GL_2(k_u)^+$ if and only if  
$\epsilon(\theta \chi^{-1},\psi_0)$ 
as well as  $\epsilon(\bar{\theta} \chi^{-1},\psi_0)$ are
independent of $\chi$; similarly for $D_u^{*+}$.

We choose a totally real number field $k$ and a place $u$ of $k$ such that the
completion of $k$ at $u$ is $k_{u}$. Let $F$ be a quadratic extension 
of $k$ which is a CM field for which $F \otimes_k k_{u} \cong L$. 
Assume that $F$ is split at all the places of $k$ 
of residue characteristic 2 except $u$ (if $u$ is of residue 
characteristic 2). Let $D$ be a quaternion division algebra 
over $k$ containing $F$ for which $D \otimes_k k_u \cong D_u$,
which  is split at all the finite places of $k$ except $u$, and 
which remains a division algebra
at all the infinite places of $k$. Such a choice of the triple $(k,F,D)$ 
is possible, as can be seen.
 
Let $\omega_{F/k} =\prod_v \omega_v$
be the quadratic character of $ {\Bbb A}_k^*/k^*$ defining $F$. 
Define
$$D^{*+}({\Bbb A}_k)= \{g \in D^*({\Bbb A}_k), g = \prod g_v| 
\omega_v(g_v) = 1, ~~\forall ~~{\rm places} ~~v ~~{\rm of} ~~k\}.$$
 Let $D^{*+}(k) = D^* \cap D^{*+}({\Bbb A}_k).$
Clearly $D^{*+}({\Bbb A}_k)$ is an open subgroup of $D^*({\Bbb A}_k)$ 
containing $D^{*+}(k)$.

Let $\pi$ be an automorphic representation of $D^*({\Bbb A}_k)$ obtained from
a character $\Theta: {\Bbb A}_F^*/F^* \rightarrow {\Bbb C}^*$ which is 
unramified at all the finite places of $F$ except $u$, and is $\theta$ 
on $F_u^*$; this is possible by lemma 3. Let ${\pi}^+ = \otimes {\pi}^+_w$ 
be an irreducible automorphic representation of $D^{*+}({\Bbb A}_k)$
contained in $\pi$  with ${\pi}^+_{u} = \pi_+$.

 Let $T$ be the
torus inside $D^*$ defined by $F$.
By lemma 1, there exists a continuous
character $\Omega = \otimes \Omega_v$ of $T(k)\backslash T({\Bbb A}_k)$ 
with $\Omega_{u}=\chi$ whose restriction to 
the center of $D({\Bbb A}_k)^*$ is the same as the central character
of $\pi$, and for which 
there exists a function $e$ in $\pi^+$ such that the integral 
$$\int_{T(k){\Bbb A}_k^*\backslash T({\Bbb A}_k)} e(t) \Omega^{-1}(t) dt,$$
is nonzero. By the theorem of Waldspurger, if $\Pi = BC(\pi)$ 
denotes the base change of 
$\pi$ to $GL_2({\Bbb A}_F)$,
$$L(\frac{1}{2}, \Pi \otimes \Omega^{-1}) \not = 0.$$
Notice that $$\Pi = BC(\pi) = Ps(\Theta \oplus \bar{\Theta}),$$
where $\bar{\Theta}$ denotes the conjugate of $\Theta$ under the 
nontrivial element of the Galois group of $F$ over $k$. 
It follows that,
$$L(s, \Pi \otimes \Omega^{-1}) =  L_F(s, \Theta \cdot \Omega^{-1})
L_F(s, \bar{\Theta} \cdot \Omega^{-1}). $$
By the condition on the central characters, one sees
 that ${\rm Ind}_F^k(\Theta \cdot \Omega^{-1})$ is a 
self-dual representation, and hence the functional equation for its
$L$-function relates itself. Therefore, since
$$ L_F(s, \Theta \cdot \Omega^{-1}) =  L_k(s, {\rm Ind}_F^k(\Theta \cdot
\Omega^{-1})),$$ 
the functional equation for $ L_F(s, \Theta \cdot \Omega^{-1}) $  relates 
itself; similarly, 
the functional equation for $ L_F(s, \bar{\Theta} \cdot \Omega^{-1}) $  
relates itself.

As $L(\frac{1}{2}, \Pi \otimes \Omega^{-1}) \not = 0$, so also are 
$L_F(\frac{1}{2}, \Theta \cdot \Omega^{-1})$ and $
L_F(\frac{1}{2}, \bar{\Theta} \cdot \Omega^{-1})$, and therefore,
in particular, the sign in their functional equations is 1:

\begin{eqnarray}
\epsilon_F(\frac{1}{2}, \Theta \cdot \Omega^{-1}) & = & 1 \\
\epsilon_F(\frac{1}{2}, \bar{\Theta} \cdot \Omega^{-1}) & =& 1.
\end{eqnarray}
The global epsilon factor is the product of the local ones:
$$\epsilon_F(\frac{1}{2}, \Theta \cdot \Omega^{-1})= 
\prod_w \epsilon_w(\frac{1}{2}, \Theta_w \cdot \Omega_w^{-1}).$$
We analyse the local epsilon factors at the various primes of $F$, and thereby
prove the theorem, in following steps.

\begin{enumerate}

\item For places $u$ of $k$ which split into  two places in $F$, say $w_1$
and $w_2$, we have
$$\epsilon_{w_1}(\frac{1}{2},\Theta_{w_1}\cdot \Omega_{w_1}^{-1})
\cdot \epsilon_{w_2}(\frac{1}{2},\Theta_{w_2}\cdot \Omega_{w_2}^{-1})
=1.$$
This follows from the following general relation:
$$\epsilon(\frac{1}{2},\chi,\psi(x))
\cdot \epsilon(\frac{1}{2},\chi^{-1},\psi(-x))
=1$$
which is what we have by the condition on the central character of
$\pi$ (which is $\Theta|_{{\Bbb A}^*_k}\cdot \omega_{F/k}$) 
being the same as the character $\Omega$ restricted to ${\Bbb A}_k^*$;
we also need to note that the additive character of ${\Bbb A}_F$ 
with respect to which we are calculating the epsilon factor is trivial
on ${\Bbb A}_k$, so is of the form $(\psi(x),\psi(-x))$ at a prime
$w$ of $k$ which splits in $F$.

\item For places $w$ of $F$ which are inert over the corresponding 
place $v$ in $k$, we note that since by
the choice of $\Omega$, the $\Omega$-period integral is nonzero, 
we have local $\Omega_w$-invariant forms on 
$\pi_w^+$ at all places. Since $\Theta$ is unramified at all finite
places away from $u$, and thus is invariant under ${\rm Gal}(F/k)$,
at such places, the representation $\pi_w$ is, up to a twist, of the
form $Ps(1,\omega)$ where $\omega$ is the quadratic character 
defining the extension $F_w$ of $k_v$. 
The next lemma now implies that
$$\epsilon_{w}(\frac{1}{2},\Theta_{w}\cdot \Omega_{w}^{-1}) =1$$
at all finite places $w$ of $F$ inert over $k$, $w \not = u$.

\item The epsilon factor,
$\epsilon_{w}(\frac{1}{2},\Theta_{w}\cdot \Omega_{w}^{-1})$ is 
 $-1$ at any  infinite place $w$ as $D$ remains a division algebra 
at such places.
 This once again follows because the period integral is nonzero, and hence
there is an invariant linear form at every infinite place, and then appealing
to the Archimedean analogue of our theorem which easily follows from
the information given in section 5.

\item Assume that $D_u = M_2(k_u)$. In this case, $D$ remains a division 
algebra at exactly the infinite primes, so they are even in number. 
We know that the
local epsilon factor at all finite places of $F$ except $u$ is 1, and is
$-1$ at the even number of infinite places. 
Since the product of the local epsilon factors at all the places of $F$
is 1, the epsilon factor 
$\epsilon_{u}(\frac{1}{2},\Theta_{u}\cdot \Omega_{u}^{-1})$
at the desired place $u$ is 1.

\item Assume $D$ is ramified at $u$ and at all the 
infinite places, which will now be odd in number.
In this case the local epsilon factor is 1 at all finite places except 
$u$, is $-1$ at all infinite places which are odd in number; 
therefore, the epsilon factor 
$\epsilon_{u}(\frac{1}{2},\Theta_{u}\cdot \Omega_{u}^{-1})=-1$ in this
case, completing the proof of the theorem.

\end{enumerate}

\begin{lemma}
Let $L$ be a quadratic extension of a local field $k_u$ of odd residue
characteristic. Then the
principal series representation $\pi = Ps(1,\omega_{L/k_u})$ splits into
two irreducible representations $\pi = \pi_+ \oplus \pi_-$ when restricted to 
$GL_2(k_u)^+ = \{x \in GL_2(k_u)| \det(x) \in NL^* \}$,
where $NL^*$ is the subgroup of $k_u^*$ of index 2 consisting of norms
from $L^*$,
 such that a character $\chi$ of $L^*$ with $\chi|_{{k_u}^*} = \omega_{L/k_u}$
appears in $\pi_+$ if and only if $\epsilon(\chi) =1$, and appears in
$\pi_-$ if and only if $\epsilon(\chi) = -1$.
\end{lemma}

\noindent{\bf Proof :} The following  was proved in lemma 3.1 of [P2]
by direct calculation in the odd residue characteristic:
$$\epsilon(\omega_{L/k_u}, \psi) \frac{{\omega_{L/k_u} \left( 
\frac{x-\bar{x}}{x_0 - \bar{x}_0} \right)}} {\left | \frac{(x-\bar{x})^2}
{x\bar{x}} \right |_{k_u}^{1/2}} = \sum_{\begin{array}{l} 
\epsilon_L(\chi, \psi_0) =  1 \\ 
\chi|_{k_u^*}\ =  \omega_{L/k_u} \end{array}} \chi(x).$$
Here $\psi$ is a nontrivial character of $k$, $x_0$ an element of
$L^*$ whose trace to $k_u$ is zero, and
 $\psi_0$ is an additive character on $L$ defined 
 by $\psi_0(x) = \psi({\rm tr}[-xx_0/2])$,
and, where, as usual, the summation (on the right) is by partial sums 
over all characters of $L^*$ of conductor $\leq n$.

It has been proved by Langlands in [L], Lemma 7.19 (the lemma with an 
`embarassing' proof),  that the character $\Theta_{\pi_+}$
of $\pi_+$ is given by  
$$\Theta_{\pi_+}(x)= \epsilon(\omega_{L/k_u}, \psi) \frac{{\omega_{L/k_u} \left ( 
\frac{x-\bar{x}}{x_0 - \bar{x}_0} \right)}} {\left | \frac{(x-\bar{x})^2}
{x\bar{x}} \right |_{k_u}^{1/2}}.$$
These two results combine to prove lemma 4.

\vspace{1 cm}

\noindent{\bf Remark :} We have stated and proved the lemma only in the odd
residue characteristic; this is enough for our purposes to prove the theorem
5 in all residue characteristics. The statement of the theorem includes
this lemma as a particular case, and therefore after we have 
proved the theorem 5 using this lemma in odd residue characteristic, 
we get a proof of lemma 4
in {\it all} residue characteristics by global method! 
Also, since the lemma of 
Langlands is available in all residue characteristics, we get 
a proof of the following lemma proved by this author only in odd residue
charateristic earlier in [P2].

\begin{lemma} Let $L$ be a quadratic extension of a local field $k_u$. 
Fix a nontrivial character $\psi$ of $k_u$, and an element $x_0 $ in $L^*$ 
whose trace to $k_u$ is zero. Define an additive character
$\psi_0$ on $L$ by $\psi_0(x) = \psi({\rm tr}[-xx_0/2])$. Then,
$$\epsilon(\omega_{L/k_u}, \psi) \frac{{\omega_{L/k_u} \left ( 
\frac{x-\bar{x}}{x_0 - \bar{x}_0} \right)}} {\left | \frac{(x-\bar{x})^2}
{x\bar{x}} \right |_{k_u}^{1/2}} = \sum_{\begin{array}{l} 
\epsilon(\chi, \psi) =  1 \\ 
\chi|_{k_u^*}\ =  \omega_{L/k_u} \end{array}} \chi(x)$$
where as usual, the summation on the right is by partial sums 
over all characters of $L^*$ of conductor $\leq n$.
\end{lemma}

\noindent{\bf Remark :} Given a quadratic extension $L$ of a local field
$k_u$, the group $L^1$ of norm 1 elements of $L^*$ can be embedded in
$SL_2(k_u)$. The two fold metaplectic cover $Mp_2(k_u)$ of $SL_2(k_u)$ splits
over $L^1$ (however, the splitting depends on the choice of a character
$\chi_0$ of $L^*$ such that $\chi_0|{k_u^*} = \omega_{L/k_u}$). 
Thus one can speak
of the restriction of the Weil representation $\omega_\psi$ 
of $Mp_2(k_u)$, associated to
a nontrivial additive character $\psi$ of $k_u$, to $L^1$. Characters $\mu$ 
of $L^1$ can be identified to characters
$\chi$ of $L^*$ with  $\chi|{k_u^*} = \omega_{L/k_u}$ by $\chi(x) = 
\mu(x/\bar{x})
\chi_0(x).$ It is a theorem of C. Moen, cf. [M],
 in odd residue characteristic, 
proved by J. Rogawski, cf. [Ro],
 in general, that a character $\chi$ of $L^*$ with 
$\chi|{k_u^*} = \omega_{L/k_u}$ appears in 
$\omega_\psi$ if and only if $\epsilon
(\chi, \psi) = 1$. Thus lemma 4 implies that the Weil representation
$\omega_\psi$ restricted to $L^1$ is closely related to a component of 
a reducible principal series. What lies behind this phenomenon- relating 
character of a linear and a non-linear group- is not
clear to this author.

\section{Symmetric and exterior squares}
The results in [P1] on the trilinear forms were refined in [P2] in the case
when $\pi_1 =\pi_2 = \pi$. In this case,
$$\pi \otimes \pi = {\rm Sym}^2(\pi) \oplus \wedge^2 (\pi).$$
Similarly for the corresponding Galois representations,
$$\sigma_\pi \otimes \sigma_\pi = {\rm Sym}^2(\sigma_\pi) \oplus \wedge^2 
(\sigma_\pi).$$

Here is the theorem about symmetric squares proved in [P2] in the odd residue
characteristic, now proved in general. It should also be pointed out 
that the methods of the Weil representation employed in [P2] could give
results about symmetric and exterior square only for $GL_2$, and not for
the representations of invertible elements of a division algebra!

\begin{theorem} Let $D_u$ be a quaternion algebra over $k_u$.
Let $\epsilon_D$ be 1 if $D_u = M_2(k_u)$, and $\epsilon_D =-1$ otherwise.
For irreducible admissible representations
$\pi$ and $\pi'$ of $D_u^*$ (which will be assumed to be infinite dimensional
if $D_u^* = GL_2(k_u)$) with $\omega_{\pi}^2\omega_{\pi'} =1$,
${\rm Sym}^2(\pi) \otimes \pi'$ has a $D_u^*$-invariant linear
form if and only if $\epsilon({\rm Sym}^2(\sigma_\pi) \otimes 
\sigma_{\pi'}) = \omega_{\pi}(-1)$, and $\epsilon({\wedge}^2
(\sigma_\pi) \otimes \sigma_{\pi'}) = \epsilon_D \omega_{\pi}(-1)$. The representation
$\wedge^2(\pi) \otimes \pi'$ has a $GL_2(k_u)$-invariant linear
form if and only if $\epsilon({\rm Sym}^2(\sigma_\pi) \otimes 
\sigma_{\pi'}) = - \omega_{\pi}(-1)$, and $\epsilon({\wedge}^2
(\sigma_\pi) \otimes \sigma_{\pi'}) = -\epsilon_D \omega_{\pi}(-1)$. 
\end{theorem}

\noindent{\bf Proof :} 
Let $k$ be a totally real 
number field with a finite place $u$ where the completion
of $k$ is $k_u$. Let $D$ be a quaternion division algebra over $k$ 
for which $D \otimes_k k_u \cong D_u$, which is split at all the other 
finite primes, and which remains a division algebra 
at all the places at infinity. Such a choice of the pair $(k,D)$ exists. 

Let $\Pi$ be an irreducible automorphic 
representation of $D^*({\Bbb A}_k)$ with $\pi$ as the local component at $u$,
unramified at all the other finite places of $k$, and certain 
representations at infinity. By lemma 1, there exists a
 representation $\Pi'$ of $D^*({\Bbb A}_k)$ with local components
$\pi'$ at $u$, functions $f_1, f_2 \in \Pi$, and $f' \in \Pi'$
such that 
$$\int_{D^* {\Bbb A}^*_k\backslash D^*({\Bbb A}_k)} f_1 f_2 f' dg \not = 0.$$
By the theorem of Harris-Kudla, $L(\frac{1}{2}, \Pi \times \Pi \times \Pi')
\not = 0.$ The proof proceeds, once again, using the factorisation 
of $L$-functions,
$$L(s, \Pi \times  \Pi \times \Pi') = L(s, {\rm Sym}^2(\Pi) \times \Pi')
L(s,\wedge^2(\Pi) \times \Pi'),$$
all of which are known to be analytic at $1/2$.
Therefore if $ L(\frac{1}{2}, \Pi \times \Pi \times \Pi') \not
= 0,$ then both 
$L(\frac{1}{2}, {\rm Sym}^2(\Pi) \times \Pi')$ and 
$ L(\frac{1}{2},\wedge^2(\Pi) \times \Pi')$ are nonzero. Since we are
once again dealing with self-dual representations, this forces
global epsilon factors to be 1:
\begin{eqnarray}
\epsilon(\frac{1}{2}, {\rm Sym}^2(\Pi) \times \Pi') & = & 1 \\ 
\epsilon(\frac{1}{2},\wedge^2(\Pi) \times \Pi') & = & 1.
\end{eqnarray}
We have the factorisation of epsilon factors,
$$\epsilon(\frac{1}{2}, {\rm Sym}^2(\Pi) \times \Pi') = 
\prod_v \epsilon(\frac{1}{2}, {\rm Sym}^2(\Pi_v) \times \Pi'_v).$$

Observe that the linear form 
$$f_1 \otimes f_2 \otimes f' \longrightarrow 
\int_{D^* {\Bbb A}^*_k\backslash D^*({\Bbb A}_k)} f_1 f_2 f' dg,$$
defines a $D^*({\Bbb A}_k)$-invariant linear form on
$\Pi \otimes  \Pi \otimes \Pi'$ which is symmetric in the first 2 variables.
By the local uniqueness of the trilinear form, the invariant form
on $\Pi_w \otimes  \Pi_w \otimes \Pi'_w$ is either symmetric or skew-symmetric
in the first 2 variables. By generalities about group representations, it
follows that the set of places $w$ of $k$ for which the invariant form
on $\Pi_w \otimes  \Pi_w \otimes \Pi'_w$ is skew-symmetric is even.
By choice, $\Pi_w$ is an unramified
principal series at all finite places $w \not = u$; in particular, we have
the condition on the central character $\omega_{\Pi_w}(-1) =1$.
By the following lemma from [P3], the local invariant forms 
on $\Pi_w \otimes  \Pi_w \otimes \Pi'_w$ at all finite
places $w \not = u$ are symmetric. 

\begin{lemma} Suppose $\pi_v$ is a principal series representation
of $GL_2(k_v)$. Then for an irreducible admissible representation
$\pi'_v$ of $GL_2(k_v)$, ${\rm Sym}^2(\pi_v) \otimes \pi'_v$ 
has a $ GL_2(k_v)$-invariant linear form if and only if $\omega_{\pi_v}(-1)=1$.

\end{lemma}

Theorem 6 now follows from the following lemma at infinity whose simple
proof will be omitted. 

\begin{lemma} Let ${\Bbb H}$ be the quaternion division algebra 
over ${\Bbb R}$.
For irreducible representations
$\pi$ and $\pi'$ of ${\Bbb H}^*$  with $\omega_{\pi}^2\omega_{\pi'} =1$,
${\rm Sym}^2(\pi) \otimes \pi'$ has a ${\Bbb H}^*$-invariant linear
form if and only if $\epsilon({\rm Sym}^2(\sigma_\pi) \otimes 
\sigma_{\pi'}) = \omega_{\pi}(-1)$, and $\epsilon({\wedge}^2
(\sigma_\pi) \otimes \sigma_{\pi'}) = - 
\omega_{\pi}(-1)$. The representation
$\wedge^2(\pi) \otimes \pi'$ has a ${\Bbb H}^*$-invariant linear
form if and only if $\epsilon({\rm Sym}^2(\sigma_\pi) \otimes 
\sigma_{\pi'}) = -\omega_{\pi}(-1)$, and $\epsilon({\wedge}^2
(\sigma_\pi) \otimes \sigma_{\pi'}) =  \omega_{\pi}(-1)$. 
\end{lemma}

\section{Archimedean Case}

The results in this paper have eventually depended on similar but much 
simpler questions for the Archimedean field (which can be assumed to be 
${\Bbb R}$). 
For this, we fix some notation and recall some standard facts without
going into the proofs of the Archimdean lemmas used 
in this paper. We will let ${\Bbb H}$
denote the quaternion division algebra over ${\Bbb R}$.

The Weil group $W_{{\Bbb C}/{\Bbb R}}$ of ${\Bbb R}$ is the normaliser of
${\Bbb C}^*$ in ${\Bbb H}^*$ and sits in the exact sequence :
$$ 0 \rightarrow {\Bbb C}^*  \rightarrow W_{{\Bbb C}/{\Bbb R}} \rightarrow
{\Bbb Z}/2{\Bbb Z} \rightarrow 0.$$

Let 
$\psi(x)= \exp(2\pi i x)$  be the character of ${\Bbb R}$, and $\psi_0(z) = 
\psi({\rm tr_ {{\Bbb C}/{\Bbb R}}  }{iz})$ be a character of ${\Bbb C}$ 
which is trivial on ${\Bbb R}$. Then for integers $n,m$, and complex number 
$s$, 
we have,

$$\begin{array}{ccccc}
\epsilon(z^n\bar{z}^m |z|^s, \psi_0) & = & (-1)^{n-m} & {\rm if} &  
n \geq m, \\
   & = & 1   & {\rm if} &  n <  m.
\end{array}$$

For $m \geq 0$, let $\sigma_m$ be the two
dimensional representation Ind$_{{\Bbb C}^*}^{W_{{\Bbb C}/{\Bbb
R}}}(z/|z|)^m$ of $W_{{\Bbb C}/{\Bbb R}}$. We have,
$$ \epsilon(\sigma_m, \psi) = 
i ^{m+1} \mbox{~~~ for ~~~~} m \geq
0.$$ 

It follows that,
\begin{eqnarray*}
\epsilon(\sigma_m \otimes \sigma_n, \psi)&  = &  
(-1) ^{m+1} \mbox{~~ for ~~} m \geq n  \\  
& = & (-1) ^{n+1} \mbox{~~~ for ~~} n \geq m.
\end{eqnarray*}

Under the Langlands correspondence, the discrete series representation $D_m$ of
$GL_2({\Bbb R})$  for $m \geq 2$, which has trivial central character 
restricted to ${\Bbb R}^{*+}$, 
corresponds to the representation 
$\sigma_{m-1}$ of the Weil group  $W_{{\Bbb C}/{\Bbb R}}$. The corresponding
representation of ${\Bbb H}^*/ {\Bbb R}^{*+} $, to be denoted by $F_{m-2}$,
 is of dimension $m-1$, and of highest weight $z^{m-2}/{|z|^{m-2}}$.

The proof of the following lemma is standard, and will therefore be omitted.
Using the information on epsilon factors given in this section, 
this is then easily seen to be equivalent to lemma 7.

\begin{lemma} Let $F_n$ denote the irreducible representation of
${\Bbb H}^*/ {\Bbb R}^{*+} $ of highest weight $n$. Then,
\begin{eqnarray*}
Sym^2(F_n) & =  & F_{2n} \oplus F_{2n-4} \oplus \cdots  \\
\Lambda^2 F_n & = & F_{2n-2} \oplus F_{2n-6} \oplus \cdots 
\end{eqnarray*}
\end{lemma}

\noindent{\bf Remark :} For the purposes of this paper, it is curious to note 
that we can assume that the groups are compact at infinity, and therefore
can deduce theorems about (discrete series representations) of noncompact
real groups from that of its compact form via the global theorems!

\section{Epilogue} 

Results in this paper were proved by using the Burger-Sarnak principle
after applying the theorem of Harris-Kudla or Waldspurger. These results 
enabled us to prove nonvanishing of central critical value of 
certain $L$-functions which played a crucial role in our
applications. 
There is a  very large literature 
on the nonvanishing results of this kind, but we do not detail them here. 
However, we take the opportunity to state the
following general questions about $L$-functions which we have answered
in this paper in some of the cases we studied here, and which might be
amenable to do by similar methods.

\vspace{2mm}

\noindent{\bf Questions:}

\begin{enumerate}

\item Given a cuspidal automorphic representation $\Pi$ of
$GL_n({\Bbb A}_k)$, and an integer $m \leq n$, is there a cuspidal
automorphic representation $\Pi'$ on $GL_m({\Bbb A}_k)$ with prescribed
discrete series behaviour at finitely many places of $k$, and with trivial
central character for $\Pi \times \Pi'$, such that 
$$L(\frac{1}{2}, \Pi \times \Pi') \not = 0.$$

\item Given a cuspidal automorphic representation $\Pi$ of
$GL_n({\Bbb A}_k)$ which is self-dual up to a twist, and an integer 
$m \leq n$, 
is there a cuspidal
automorphic representation $\Pi'$ on $GL_m({\Bbb A}_k)$ which 
is self-dual up to a twist, has prescribed discrete series 
behaviour at finitely many places of $k$, such that 
$\Pi \times \Pi'$ is self-dual, and  such that
$$L(\frac{1}{2}, \Pi \times \Pi') \not = 0.$$
\end{enumerate}

\vspace{2mm}

\noindent{\bf Acknowledgement :} The author thanks T.N. Venkataramana
for discussions on the Burger-Sarnak `principle' on several occasions, 
and also thanks U.K. Anandavardhanan for some helpful discussions.


\vspace{1cm}

\noindent School of Mathematics, Tata Institute of Fundamental Research,

\noindent Colaba, Mumbai-400005, INDIA

\noindent Email: dprasad@math.tifr.res.in

\end{document}